\documentclass[11pt,reqno]{amsart}
\usepackage{indentfirst, amssymb,amsmath,amsthm}
\usepackage[colorlinks]{hyperref}
\newtheorem{theo}{Theorem}[section]
\newtheorem{exm}{Example}[section]
\newtheorem{lem}{Lemma}
\newtheorem*{prf}{Proof}
\newtheorem{pro}{Proposition}[section]
\newtheorem{cor}{Corollary}[section]

\newtheorem{res}{Result}[section]
\newtheorem{defi}{Definition}
\newtheorem{obs}{Observation}
\newtheorem{note}{Note}

\newcommand{\be}{\begin{equation}}
\newcommand{\ee}{\end{equation}}
\newcommand{\beas}{\begin{eqnarray*}}
\newcommand{\eeas}{\end{eqnarray*}}
\newcommand{\bea}{\begin{eqnarray}}
\newcommand{\eea}{\end{eqnarray}}

\numberwithin{equation}{section}

\begin{document}

\setlength{\unitlength}{1mm} \baselineskip .52cm
%\large
\setcounter{page}{1}
\pagenumbering{arabic}
\title[Some Properties of Lebesgue Fuzzy Metric Spaces]{Some Properties of Lebesgue Fuzzy Metric Spaces}

\author[Sugata Adhya and A. Deb Ray]{Sugata Adhya and A. Deb Ray}

\address{Department of Mathematics, The Bhawanipur Education Society College. 5, Lala Lajpat Rai Sarani, Kolkata 700020, West Bengal, India.}
\email {sugataadhya@yahoo.com}

\address{Department of Pure Mathematics, University of Calcutta. 35, Ballygunge Circular Road, Kolkata 700019, West Bengal, India.}
\email {debrayatasi@gmail.com}

\maketitle

\begin{abstract}
In this paper, we establish a sequential characterisation of Lebesgue fuzzy metric and explore the relationship between Lebesgue, weak $G$-complete and compact fuzzy metric spaces. We also discuss the Lebesgue property of several well-known fuzzy metric spaces.\\
\end{abstract}

\noindent{\textit{AMS Subject Classification:} 54A40, 54E35, 54E40.}\\
{\textit{Keywords:} {Fuzzy metric space, Lebesgue property, weak $G$-complete.}} 

\section{\textbf{Introduction}}

The theory of fuzzy metric spaces, proposed by George and Veeramani \cite{ver1}, is widely accepted as a consistent notion for metric fuzziness. It is a slight modification of the one due to Kramosil and Michalek \cite{kmf}. Throughout the paper, this is the only notion of fuzzy metric we will be working on. It should be noted that every fuzzy metric gives rise to a metrizable topology that allowed the researchers to adopt several concepts from metric spaces in this fuzzy setting. In particular, Gregori, Romaguera, and Sapena \cite{com2} introduced a notion similar to the Lebesgue number in the realm of fuzzy metric spaces. 

In the theory of metric spaces, the Lebesgue number lemma states that every open cover $\mathcal U$ of a compact metric space $(X, d)$ corresponds to a positive number $\delta$ such that any subset of $X$ having diameter less than $\delta$ gets contained in some member of $\mathcal U.$ This $\delta$ is called a Lebesgue number for $\mathcal U.$ The property of having such positive real numbers for every open cover is called the Lebesgue property for metric spaces. It is important to note that one can find non-compact metric spaces (e.g. consider the set of positive integers endowed with discrete topology) that satisfy Lebesgue property. In fact, the study of metric spaces having Lebesgue property (precisely, Lebesgue metric spaces) is an interesting problem in the theory of metric spaces. For details one may refer to \cite{gb} and references therein. 

In 2001, Gregori, Romaguera and Sapena \cite{com2} gave a satisfactory extension to the notion of Lebesgue property for fuzzy metric spaces and characterized it in terms of uniform continuity, equinormality and uniformity. They ensured the existence of a non-standard Lebesgue fuzzy metric that made Lebesgue property worth studying in the realm of fuzzy metric setting. Unfortunately, Lebesgue fuzzy metric spaces didn't get much attention of researchers, later on. Recently in \cite{2}, we discussed several new characterizations for Lebesgue fuzzy metric spaces and this paper is a continuation of that. 

In this paper, we provide a sequential characterization for Lebesgue fuzzy metric and employ it to study the Lebesgue property of some well-known fuzzy metric spaces. In what follows, we use the characterization to study the relationship between Lebesgue, weak $G$-complete and compact fuzzy metric spaces. 

Throughout the paper, $\mathbb R$ and $\mathbb N$ will stand for the sets of real numbers and positive integers, respectively.

\section{\textbf{Preliminaries}}

In this section, we recall a series of definitions and some related results on fuzzy metric spaces that will be required subsequently. For undefined terms related to general topology, we refer to \cite{will}.

\begin{defi}
\normalfont \cite{tn} Let $*$ be a binary operation on $I=[0,1]$ which is associative, commutative and continuous on $I\times I$. Then $*$ is said to be a continuous $t$-norm, if 

a) $\forall~a\in[0,1],~a*1=a;$

b) $\forall~a,b,c,d\in[0,1],~a\le b,~c\le d\implies a*c\le b*d.$
\end{defi}

\begin{defi}
\normalfont \cite{ver1,ver2} Given a non-empty set $X,$ a continuous $t$-norm $*$ and a mapping $M:X\times X\times(0,\infty)\to[0,1],$ the ordered pair $(M,*)$ is said to be a fuzzy metric on $X$ if, for all $x,y\in X$ and $t>0,$ the following conditions hold:

a) $M(x,y,t)>0;$

b) $M(x,y,t)=1\iff x=y;$

c) $M(x,y,t)=M(y,x,t);$

d) $M(x,y,t)*M(y,z,s)\le M(x,z,t+s);$

e) $M(x,y,.):(0,\infty)\to[0,1]$ is continuous.

In this case, $(X, M,*)$ is said to be a fuzzy metric space.
\end{defi}

It is easy to note from the above axioms that given two elements $x,y$ in a fuzzy metric space $(X,M,*),$ the mapping $t\mapsto M(x,y,t)$ is increasing on $(0,\infty).$

\begin{res}
\normalfont\cite{ver1} Let $(X,M,*)$ be a fuzzy metric space. Then $\{B_M(x,r,t):x\in X,r\in(0,1),t>0\},$ where $B_M(x,r,t)=\{y\in X:M(x,y,t)>1-r\},$ forms a base for some topology $\tau_M$ on $X.$
\end{res}

\begin{defi}
\normalfont $\tau_M$ is called the topology induced by $(M,*).$
\end{defi}

\begin{defi}
\normalfont \cite{ver1} Let $(X,d)$ be a metric space. If $M_d:X\times X\times(0,\infty)\to[0,1]$ is defined for all $x,y\in X$ and $t>0$ by $$M_d(x,y,t)=\frac{t}{t+d(x,y)}$$ then $(M_d,~\cdot),$ `$\cdot$' being the usual multiplication on $[0,1]$, defines a fuzzy metric on $X.$ It is called the standard fuzzy metric induced by $d.$
\end{defi}

\begin{res}\label{res1}
\normalfont \cite{ver1} If $(X,d)$ is a metric space, then $\tau_{M_d}=\tau(d),$ where $\tau(d)$ denotes the topology induced by the metric $d.$
\end{res}

\begin{defi}
\normalfont \cite{stat} A fuzzy metric space $(X,M,*)$ is said to be stationary if for all $x,y\in X,~t\mapsto M(x,y,t)$ defines a constant mapping on $(0,\infty).$
\end{defi}

George and Veeramani \cite{ver1} initiated the study of convergence of sequences for fuzzy metric spaces. A sequence $(x_n)$ in a fuzzy metric space $(X,M,*)$ converges to $x$ (\textit{resp.} clusters), if it does so in $(X,\tau_M).$ 

\begin{theo}
\normalfont \cite{ver1} A sequence $(x_n)$ in a fuzzy metric space $(X,M,*)$ converges to $x\in X$ if and only if $\lim\limits_{n\to\infty}M(x_n,x,t)= 1,~\forall~t>0.$.
\end{theo}

\begin{defi}
\normalfont \cite{ver1} A sequence $(x_n)$ in a fuzzy metric space $(X,M,*)$ is said to be Cauchy if for $\epsilon\in(0,1)$ and $t>0,$ there exists $k\in\mathbb N$ such that  $M(x_m,x_n,t)>1-\epsilon,~\forall~m,n\ge k.$

A fuzzy metric space, in which every Cauchy sequence converges, is said to be complete.
\end{defi}

\begin{defi}
\normalfont \cite{uni} A fuzzy metric space $(X,M,*)$ is said to be precompact if for $r\in(0,1)$ and $t>0,$ there exists a finite subset $A$ of $X$ such that $X=\bigcup\limits_{x\in A}{B_M(x,r,t)}.$
\end{defi}

\begin{pro}\label{prec}
\normalfont \cite{pc} A metric space $(X,d)$ is precompact if and only if the standard fuzzy metric space $(X,M_d,\cdot)$ is precompact.
\end{pro}

\begin{lem}\label{uni1}
\normalfont \cite{uni} A fuzzy metric space $(X,M,*)$ is precompact if and only if every sequence in $X$ has a Cauchy subsequence.
\end{lem}

In \cite{uni}, Gregori and Romaguera introduced compactness  for fuzzy metric spaces in the most obvious way: A fuzzy metric space $(X,M,*)$ is compact if so is $(X,\tau_M)$ as a topological space. They characterized compact fuzzy metric spaces as follows:

\begin{theo}\label{uni2}
\normalfont \cite{uni} A fuzzy metric space $(X,M,*)$ is compact if and only if it is precompact and complete.
\end{theo}

\section{\textbf{Sequential Characterization for Lebesgue Property}}

\begin{defi}
\normalfont \cite{com2} A fuzzy metric space $(X,M,*)$ is said to have the Lebesgue property if given an open cover $\mathcal G$ of $(X,\tau_M),$ there exist $r\in(0,1),~t>0$ such that $\{B_M(x,r,t):x\in X\}$ refines $\mathcal G.$ We call such fuzzy metric spaces \textit{Lebesgue}.
\end{defi}

\begin{pro}\label{pro1}
\normalfont \cite{com2} Let $(X,d)$ be a metric space. Then $(X,d)$ is Lebesgue if and only if $(X,M_d,\cdot)$ is Lebesgue.
\end{pro}

\begin{defi}
\normalfont \cite{com2} A fuzzy metric space $(X,M,*)$ is said to be equinormal if for given nonempty, closed subsets $B$ and $C$ of $(X,\tau_M)$ with $B\cap C=\emptyset,$ there exists $s>0$ such that $\sup\{M(b,c,s):b\in B,c\in C\}<1.$
\end{defi}

Several characterizations of the Lebesgue property for fuzzy metric spaces have been discussed in \cite{2} and \cite{com2}. In particular, it has been shown in \cite{com2} that a fuzzy metric space is Lebesgue if and only if it is equinormal. In what follows, we give a sequential characterization for Lebesgue fuzzy metric spaces. 

To attain the requirement of our main result, we first extend the notion of pseudo-Cauchy sequences in fuzzy metric setting. 

Recall that, a sequence $(x_n)$ in a metric space $(X,d)$ is pseudo-Cauchy if given $\epsilon>0$ and $k\in\mathbb N,$ there exist $j,n~(>k)\in\mathbb N$ with $j\ne n$ such that $d(x_j,x_n)<\epsilon.$ We propose the notion of fuzzy pseudo-Cauchy sequence as follows:

\begin{defi}
\normalfont A sequence $(x_n)$ in a fuzzy metric space $(X,M,*)$ is said to be fuzzy pseudo-Cauchy if given $\epsilon\in(0,1),~t>0$ and $k\in\mathbb N,$ there exist $j,~n~(>k)\in\mathbb N$ with $j\ne n$ such that $M(x_j,x_n,t)>1-\epsilon.$
\end{defi}

Clearly, a Cauchy sequence in a (fuzzy) metric space is (fuzzy) pseudo-Cauchy, however, the converse may fail.

\begin{pro}\label{th1}
\normalfont Let $(X,d)$ be a metric space. A sequence $(x_n)$ in $(X,M_d,\cdot)$ is fuzzy pseudo-Cauchy if and only if $(x_n)$ is pseudo-Cauchy in $(X,d).$
\end{pro}

\begin{prf}
\normalfont Consider a fuzzy pseudo-Cauchy sequence $(x_n)$ in $(X,M_d,\cdot).$

Choose $\epsilon\in(0,1)$ and $k\in\mathbb N.$ Then there exist $j,n~(>k)\in\mathbb N$ with $j\ne n$ such that $M_d(x_j,x_n,1-\epsilon)>1-\epsilon.$

i.e., $d(x_j,x_n)<\epsilon.$ Thus $(x_n)$ is pseudo-Cauchy in $(X,d).$

\textit{Conversely,} let $(x_n)$ be a pseudo-Cauchy sequence in $(X,d).$

Choose $\epsilon\in(0,1),~t>0$ and $k\in\mathbb N.$ Then there exist $j,n~(>k)\in\mathbb N$ with $j\ne n$ such that $d(x_j,x_n)<\frac{t\epsilon}{1-\epsilon}.$

Consequently, $\frac{t}{t+d(x_j,x_n)}>1-\epsilon\implies M_d(x_j,x_n,t)>1-\epsilon.$ Thus $(x_n)$ is fuzzy pseudo-Cauchy in $(X,M_d,\cdot).\qed$ 
\end{prf}

\begin{exm}\label{exm1}
\normalfont Consider the non-standard fuzzy metric space $(\mathbb N,M,*)$ \cite{com2} where $a*b=ab,~\forall~a,b\in[0,1]$ and for $x,y\in\mathbb N,~t>0,$ $$M(x,y,t)=\begin{cases}1 &\text{ if }x=y\\\frac{1}{xy} &\text{ otherwise}\end{cases}.$$

Then $(1,~2,~1,~3,~1,~4,\cdots)$ is a fuzzy pseudo-Cauchy sequence in $(\mathbb N,M,*)$ which is not Cauchy.
\end{exm}

We are now at a stage to discuss the main result of this section.

\begin{theo}\label{thmn}
\normalfont Let $(X,M,*)$ be a fuzzy metric space. Then $(X,M,*)$ is Lebesgue if and only if every fuzzy pseudo-Cauchy sequence in $(X,M,*)$ having distinct terms has a cluster point in $(X,\tau_M).$
\end{theo}

\begin{prf}
\normalfont Let $(X,M,*)$ be Lebesgue. 

Choose a fuzzy pseudo-Cauchy sequence $(x_n)$ having distinct terms in $X.$ Then there exists a strictly increasing sequence $(r_n)$ of natural numbers such that $$M\left(x_{r_{2n-1}},x_{r_{2n}},\frac{1}{n+1}\right)>1-\frac{1}{n+1},~\forall~n\in\mathbb N.$$

If possible, let none of $(x_{r_{2n-1}})$ and $(x_{r_{2n}})$ has cluster point in $(X,\tau_M)$. Then, $B=\{x_{r_{2n-1}}:n\in\mathbb N\}$ and $C=\{x_{r_{2n}}:n\in\mathbb N\}$ are disjoint, closed subsets of $(X,\tau_M).$

Since $(X,M,*)$ is equinormal, being Lebesgue, there exists $s>0$ such that $$\sup\{M(b,c,s):b\in B, c\in C\}=p\cdots(*)$$ where $p<1.$

Choose, $k\in\mathbb N$ such that $\frac{1}{k}<\min\{s,1-p\}.$

Then $M(x_{r_{2n-1}},x_{r_{2n}},s)\ge M(x_{r_{2n-1}},x_{r_{2n}},\frac{1}{n+1})>1-\frac{1}{n+1}>p,~\forall~n\ge k,$ a contradiction to $(*).$

Thus, at least one of $(x_{r_{2n-1}})$ or $(x_{r_{2n}})$ has a cluster point in $(X,\tau_m)$, which establishes the fact that $(x_n)$ has a cluster point in $(X,\tau_m).$

\textit{Conversely,} let the condition hold. If possible, let $(X,M,*)$ be not a Lebesgue fuzzy metric space. Then there exists an open cover $\mathcal{G}=\{U_\lambda:\lambda\in\Lambda\}$ of $(X,\tau_M)$ such that for no $r\in(0,1)$ and $s>0,$ $\{B(x,r,s):x\in X\}$ refines $\mathcal{G}.$ Thus for each $n\ge1,$ there exists $x_{2n-1}\in X$ such that $B\left(x_{2n-1},\frac{1}{n+1},\frac{1}{n+1}\right)\not\subset U_\lambda,~\forall~\lambda\in\Lambda.$ Since $\mathcal G$ covers $(X,\tau_M),~x_{2n-1}\in U_{\lambda_n}$ for some $\lambda_n\in\Lambda.$ Thus $\exists$ $x_{2n}\in B\left(x_{2n-1},\frac{1}{n+1},\frac{1}{n+1}\right)\backslash U_{\lambda_{n}}$ and consequently, $M(x_{2n-1},x_{2n},\frac{1}{n+1})>1-\frac{1}{n+1},~\forall~n\ge1\cdots(1).$

We first show that $(x_n)$ is a fuzzy pseudo-Cauchy sequence. Let $\epsilon\in(0,1),t>0$ and $k>1.$ Choose $q>k$ such that $\frac{1}{q}<\min\{\epsilon,t\}.$ Then $2q-1,2q>k$ and $M(x_{2q-1},x_{2q},t)\ge M\left(x_{2q-1},x_{2q},\frac{1}{q+1}\right)>1-\frac{1}{q+1}>1-\epsilon.$ Thus $(x_n)$ is fuzzy pseudo-Cauchy.

We now show that $(x_n)$ has a fuzzy pseudo-Cauchy subsequence $(x_{r_n})$ of distinct terms. 

\textit{Case I:} Suppose $(x_n)$ does not have any constant subsequence. We proceed by induction. 

Set $x_{r_1}=x_1$ and $x_{r_2}=x_2.$ For chosen $\{x_{r_1},x_{r_2},x_{r_3},x_{r_4},\cdots,x_{r_{2k-1}},x_{r_{2k}}\}$ find $p>r_{2k}$ such that $x_{2p-1},x_{2p}\notin\{x_{r_1},x_{r_2},x_{r_3},x_{r_4},\cdots,x_{r_{2k-1}},x_{r_{2k}}\}$ and set $x_{r_{2k+1}}=x_{2p-1},~x_{r_{2k+2}}=x_{2p}.$ Thus we obtain a subsequence $(x_{r_n})$ of $(x_n)$ having distinct terms.

Choose $t>0,~\epsilon\in(0,1)$ and $k\in\mathbb N.$ Find $q>k$ such that $\frac{1}{q}<\min\{\epsilon,t\}.$ Then $M(x_{2{r_q}-1},x_{2{r_q}},t)\ge M(x_{2{r_q}-1},x_{2{r_q}},\frac{1}{q+1})\ge M(x_{2{r_q}-1},x_{2{r_q}},\frac{1}{r_q+1})>1-\frac{1}{r_q+1}\ge1-\frac{1}{q+1}>1-\epsilon.$ Consequently $(x_{r_n})$ is fuzzy pseudo-Cauchy.

\textit{Case II:} Suppose $(x_n)$ has a constant subsequence $(x_{r_n}),$ where $x_{r_n}=a,~\forall~n\ge1.$ 

By setting $$m_n=\begin{cases}r_n-1,&\text{if $r_n$ is even}\\r_n+1,&\text{if $r_n$ is odd}\end{cases}$$ we see that for chosen $n\in\mathbb N,~\exists~k\in\mathbb N$ such that $\{r_n,m_n\}=\{2k-1,2k\}.$ Since $(r_n)$ defines a strictly increasing sequence of natural numbers, so does $(m_n).$ Thus $(x_{m_n})$ forms a subsequence of $(x_n).$

We first show that, $(x_{m_{2n}})$ is fuzzy pseudo-Cauchy. 

Using equation $(1)$ we see that, $\forall~n\in\mathbb N,~M\left(x_{r_{2n}},x_{m_{2n}},\frac{1}{n+1}\right)>1-\frac{1}{n+1},$ that is, $M(a,x_{m_{2n}},\frac{1}{n+1})>1-\frac{1}{n+1}.$ Thus $(x_{m_{2n}})$ is convergent and hence, is fuzzy pseudo-Cauchy.

Since $\lim\limits_{n\to\infty}x_{m_{2n}}=a,$ and $x_{m_{2n}}\ne x_{r_{2n}},~\forall~n\in\mathbb N,$ it follows that $(x_{m_{2n}})$ has no constant subsequence. Thus, in view of {Case I}, it must have a fuzzy pseudo-Cauchy subsequence of distinct terms. 

Consequently, in any case, $(x_n)$ has a fuzzy pseudo-Cauchy subsequence of distinct terms.

Thus, in view of the hypothesis, $(x_n)$ must have a cluster point $z$ in $(X,\tau_M).$

Clearly $z\in U_\lambda$ for some $\lambda\in\Lambda.$ Since $U_\lambda$ is open, there exists $r\in(0,1),s>0$ such that $B_M(z,r,s)\subset U_\lambda.$

Since $*$ is continuous, there exists $r'\in(0,1)$ with $r'<r$ such that $(1-r')*(1-r')*(1-r')>1-r.$ 

Also $z$ being a cluster point of $(x_n),$ there is a natural number $p$ satisfying $\frac{1}{p}<\min\{r',\frac{s}{3}\}$ such that at least one of $x_{2p}$ and $x_{2p-1}$ belongs to $B_M(z,r',\frac{s}{3}).$ Set $y$ to be a point among $x_{2p}$ and $x_{2p-1}$ such that it lies in $B_M(z,r',\frac{s}{3}).$

Note for $w\in B_M(x_{2p-1},\frac{1}{p+1},\frac{1}{p+1}),$ we have $M(w,z,s)\ge M(w,x_{2p-1},\frac{s}{3})*M(x_{2p-1},y,\frac{s}{3})*M(y,z,\frac{s}{3})\ge M(w,x_{2p-1},\frac{1}{p+1})*M(x_{2p-1},y,\frac{1}{p+1})*M(y,z,\frac{s}{3})\ge(1-\frac{1}{p+1})*(1-\frac{1}{p+1})*(1-r')\ge(1-r')*(1-r')*(1-r')>1-r,$ that is, $w\in B_M(z,r,s).$

Thus $B_M(x_{2p-1},\frac{1}{p+1},\frac{1}{p+1})\subset B_M(z,r,s)\subset U_\lambda,$ a contradiction.

So $(X,M,*)$ is Lebesgue$.\qed$
\end{prf}

\begin{exm}\label{exm2}
\normalfont It is worth recalling, at this stage, that $(\mathbb N,M,*),$ defined in Example \ref{exm1}, forms a non-standard, Lebesgue fuzzy metric space \cite{com2}. 

In fact, Theorem \ref{thmn} can be employed to realize that $(\mathbb N,M,*)$ is Lebesgue: Choose $\epsilon=\frac{1}{2}$ and $t>0.$ Then for no $x,y~(x\ne y)\in\mathbb N$ we can have $M(x,y,t)>1-\epsilon.$ So, there is no fuzzy pseudo-Cauchy sequence in $(\mathbb N,M,*).$ Consequently, $(\mathbb N,M,*)$ is Lebesgue.
\end{exm}

Before proceeding further, we note from pseudo-Cauchy characterization of the Lebesgue property that the class of Lebesgue fuzzy metric spaces resides strictly in-between the classes of compact and complete fuzzy metric spaces.

\begin{exm}\label{exnl1}
\normalfont For $X=(0,\infty),$ define $M:X^2\times(0,\infty)\to[0,1]$ by $M(x,y,t)=\frac{\min\{x,y\}}{\max\{x,y\}},~\forall~x,y\in X,~t>0.$ It has been shown in \cite{uc} that, $(X,M,\cdot)$ forms a complete fuzzy metric space which is not compact.

We now show that, $(X,M,\cdot)$ is \textit{not even} Lebesgue.

Set $a_n=n,~\forall~n\in\mathbb N.$ Then $\lim\limits_{n\to\infty}(a_n,a_{n+1},t)=\lim\limits_{n\to\infty}\frac{n}{n+1}=1,~\forall~t>0.$ So, given $\epsilon\in(0,1),t>0$ and $k\in\mathbb N,~\exists~p,~q~(p\ne q)>n$ such that $M(x_p,x_q,t)>1-\epsilon.$ Thus, $(a_n)$ is a fuzzy pseudo-Cauchy sequence of distinct terms in $X.$

If possible, let $c$ be a cluster point of $(a_n).$ Then there exists a subsequence $(a_{r_n})$ of $(a_n)$ that converges to $c$ with respect to the topology $\tau_M.$ Note that, $\exists~k\in\mathbb N$ such that $a_{r_n}>c,~\forall~n\ge k,$ whence $\lim\limits_{n\to\infty}M(a_{r_n},c,t)=\lim\limits_{n\to\infty}\frac{c}{r_n}=0,~\forall~t>0.$ Thus $(a_{r_n})$ cannot converge to $c,$ a contradiction. Consequently, $(a_n)$ has no cluster point.

So, in view of Theorem \ref{thmn}, $(X,M,\cdot)$ is not Lebesgue.
\end{exm}

\begin{exm}
\normalfont For $X=[0,\infty),$ define $M:X^2\times(0,\infty)\to[0,1]$ by $M(x,y,t)=\frac{\min\{x,y\}+t}{\max\{x,y\}+t},~\forall~x,y\in X,~t>0.$ Then $(X,M,\cdot)$ forms a complete fuzzy metric space \cite{uc}. Arguing as Example \ref{exnl1}, it can be shown that $(X,M,\cdot)$ is not Lebesgue.
\end{exm}

\begin{exm}\label{exl2}
\normalfont Let $\phi:(0,\infty)\to(0,1]$ be a function such that $\phi(t)=t,~t\le1$ and $\phi(t)=1,\text{otherwise}.$ For $X=(0,\infty),$ define $M:X^2\times(0,\infty)\to[0,1]$ by $$M(x,y,t)=\begin{cases}1&x=y,\\\frac{\min\{x,y\}}{\max\{x,y\}}.\phi(t)&x\ne y.\end{cases}$$ It has been shown in \cite{com} that, $(X,M,\cdot)$ forms a complete fuzzy metric space.

We now show that $(X,M,\cdot)$ is, \textit{in fact}, Lebesgue.

Choose a fuzzy pseudo-Cauchy sequence $(a_n)$ of distinct terms in $(X,M,\cdot)$. Then there exists a strictly increasing sequence $(r_n)$ of natural numbers such that $$M\left(a_{r_{2n-1}},a_{r_{2n}},\frac{1}{n+1}\right)>1-\frac{1}{n+1},~\forall~n\in\mathbb N.$$

Note for chosen $t>0$ and $\epsilon\in(0,1),$ we can find $p\in\mathbb N$ such that $\frac{1}{p}<\min\{\epsilon,~t\}.$ Then $M(a_{r_{2n-1}},a_{r_{2n}},t)\ge M(a_{r_{2n-1}},a_{r_{2n}},\frac{1}{n+1})>1-\frac{1}{n+1}>1-\epsilon,~\forall~n\ge p,$ whence $\lim\limits_{n\to\infty} M(a_{r_{2n-1}},a_{r_{2n}},t)=1,~\forall~t>0.$

However $\lim\limits_{n\to\infty}M(a_{r_{2n-1}},a_{r_{2n}},\frac{1}{2})=\frac{1}{2}\times\lim\limits_{n\to\infty}\frac{\min\{a_{r_{2n-1}},a_{r_{2n}}\}}{\max\{a_{r_{2n-1}},a_{r_{2n}}\}}=\frac{1}{2}\times\lim\limits_{n\to\infty}M(a_{r_{2n-1}},a_{r_{2n}},1),$ a contradiction.

Thus no such fuzzy pseudo-Cauchy sequence $(a_n)$ exist in $(X,M,\cdot).$

So, in view of Theorem \ref{thmn}, $(X,M,\cdot)$ is Lebesgue.
\end{exm}

\begin{exm}
\normalfont For $X=(0,1),$ define $M:X^2\times(0,\infty)\to[0,1]$ by $$M(x,y,t)=\begin{cases}1&x=y,\\xy.\phi(t)&x\ne y\end{cases}$$ where $\phi$ is defined in Example \ref{exl2}.

It has been shown in \cite{com} that $(X,M,\cdot)$ forms a complete fuzzy metric space . Arguing as Example \ref{exl2}, we can see that $(X,M,\cdot)$ is, \textit{in fact}, Lebesgue.
\end{exm}

\section{\textbf{Weak $G$-Completeness \textit{versus} Lebesgue Property}}

In this section, we investigate the relationship between weak $G$-completeness and Lebesgue property for (fuzzy) metric spaces. We start by recalling the following weaker notion than Cauchy sequences, due to M. Grabiec \cite{ogc}.

\begin{defi}
\normalfont A sequence $(x_n)$ in a fuzzy metric space $(X,M,*)$ is said to be $G$-Cauchy if for each $t>0$ and $p\in\mathbb N,~\lim\limits_{n\to\infty}M(x_{n},x_{n+p},t)=1,$ \textit{or equivalently,} $\lim\limits_{n\to\infty}M(x_{n},x_{n+1},t)=1,~\forall~t>0.$
\end{defi}

Tirado, in \cite{tir}, proposed the notion of $G$-Cauchyness for metric spaces:

\begin{defi}
\normalfont A sequence $(x_n)$ in a metric space $(X,d)$ is said to be $G$-Cauchy if for each $p\in\mathbb N,~\lim\limits_{n\to\infty}d(x_{n},x_{n+p})=0,$ \textit{or equivalently,} $\lim\limits_{n\to\infty}d(x_{n},x_{n+1})=0,~\forall~t>0.$
\end{defi}

\begin{defi}
\normalfont \cite{gc} A (fuzzy) metric space $X$ is said to be 

i) weak $G$-complete if every $G$-Cauchy sequence in $X$ has a cluster point in it;

ii) $G$-complete if every $G$-Cauchy sequence in $X$ converges in it.
\end{defi}

Clearly, $$G\text{-completeness }\rightarrow\text{ weak }G\text{-completeness }\rightarrow\text{completeness}$$ though the implications cannot be reversed as is shown in \cite{gc}.

\begin{note}
\normalfont $(X,M,\cdot),$ in Example \ref{exnl1}, is not weak $G$-complete, since $(n)$ is a $G$-Cauchy sequence in $X$ without any cluster point. 
\end{note}

\begin{obs}
\normalfont It is known that $(X,M,\cdot),$ where $X=[0,1]$ and $M(x,y,t)=\frac{\min\{x,y\}+t}{\min\{x,y\}+t},$ $\forall~x,y\in X,t>0,$ defines a compact, non-$G$-complete fuzzy metric space \cite{gc}. Thus a Lebesgue fuzzy metric space may not be $G$-complete.
\end{obs}

Let us recall the following results before proceeding further:

\begin{pro}\label{pro2}
\normalfont \cite{gc} Let $(X,d)$ be a metric space. Then $(X,d)$ is weak $G$-complete if and only if the standard fuzzy metric space $(X,M_d,\cdot)$ is weak $G$-complete.
\end{pro}

It is observed in \cite{gc} that every compact metric space is weak $G$-complete. A stronger result can be realized from the succeeding discussion.

\begin{theo}\label{thgc}
\normalfont A Lebesgue fuzzy metric space is weak $G$-complete.
\end{theo}

\begin{prf}
\normalfont Let $(X,M,*)$ be a Lebesgue fuzzy metric space and $(x_n)$ be a $G$-Cauchy sequence in $X.$

If $(x_n)$ has a constant subsequence, then it must have a cluster point in $X.$

So we assume that, $(x_n)$ has no constant subsequence. 

We proceed by induction. Choose $x_{r_1},x_{r_2}$  from the sequence such that $r_1<r_2,~x_{r_1}\ne x_{r_2}$ and $M(x_{r_1},x_{r_2},\frac{1}{2})>1-\frac{1}{2}.$

Next, for chosen $x_{r_1},x_{r_{2}},\cdots,x_{r_{2k-1}},x_{r_{2k}},$ find $x_{r_{2k+1}},x_{r_{2k+2}}\notin\{x_{r_1},x_{r_{2}},\cdots,x_{r_{2k-1}},x_{r_{2k}}\}$ such that ${r_{2k}}<{r_{2k+1}}<{r_{2k+2}},$  $x_{r_{2k+1}}\ne x_{r_{2k+2}}$ and $M(x_{r_{2k+1}},x_{r_{2k+2}},$ $\frac{1}{2k+2})>1-\frac{1}{2k+2}.$

Clearly $(x_{r_n})$ defines a subsequence of distinct terms.

Choose $\epsilon\in(0,1)$ and $t>0.$ Then for any $k\in\mathbb N$ satisfying $\frac{1}{k}<\min\{\epsilon,t\},$ we have $M(x_{r_{2k+1}},x_{r_{2k+2}},t)\ge M(x_{r_{2k+1}},x_{r_{2k+2}},\frac{1}{k})\ge M(x_{r_{2k+1}},x_{r_{2k+2}},\frac{1}{2k+2})>1-\frac{1}{2k+2}>1-\epsilon.$

Consequently, $(x_{r_n})$ is a fuzzy pseudo-Cauchy sequence. So by hypothesis, $(x_{r_n}),$ and hence $(x_n),$ has a cluster point in $X.$ Hence the result follows$.\qed$
\end{prf}

The following corollary is immediate from Proposition \ref{pro1} and Proposition \ref{pro2}.

\begin{cor}
\normalfont A Lebesgue metric space is weak $G$-complete.
\end{cor}

In view of Theorem \ref{uni2}, it is now clear that the class of Lebesgue fuzzy metric spaces $\mathcal{L}$ lies in-between the classes of compact fuzzy metric spaces $\mathcal{K}$ and weak $G$-complete fuzzy metric spaces $\mathcal{G}$. In what follows, we show that $\mathcal{K}\subsetneq\mathcal{L}\subsetneq\mathcal{G}.$

\begin{exm}\label{ex27}
\normalfont (A weak $G$-complete, non-Lebesgue metric space) Let $X=\{n:n\in\mathbb N\}\cup\left\{n+\frac{1}{n}:n\in\mathbb N\right\}$ and $d$ be the usual metric on $\mathbb R$ restricted to $X\times X.$ 

\textit{$(X,d)$ is not Lebesgue:} Clearly, $\tau_d$ is the discrete topology on $X.$ Thus $\left\{\{x\}:x\in X\right\}$ is an open cover of $X$ without any Lebesgue number. Consequently $(X,d)$ is not Lebesgue.

\textit{$(X,d)$ is weak $G$-complete:} It suffices to show that the only sequences which are $G$-Cauchy are those that contain a constant subsequence. 

If possible, let there exist a $G$-Cauchy sequence $(x_n)$ in $X$ which does not have a constant subsequence. Then there exists a subsequence $(x_{r_n})$ of $(x_n)$ having distinct terms such that $d(x_{r_{n+1}},x_{r_n})<\frac{1}{3},~\forall~n\in\mathbb N,$ a contradiction. Hence $(X,d)$ is weak $G$-complete.
\end{exm}

In view of the last example, the following observation is immediate from Proposition \ref{pro1} and Proposition \ref{pro2}:

\begin{obs}
\normalfont $\mathcal{L}\subsetneq\mathcal{G}.$
\end{obs}

\begin{exm}
\normalfont (A non-compact, Lebesgue fuzzy metric space) Let $X=\left\{\frac{1}{2^n}:n\ge2\right\}\cup\left[\frac{1}{2},1\right].$ It has been shown in \cite{gc} that the stationary fuzzy metric space $(X,M,\cdot),$ where $M(x,y,t)=\frac{\min\{x,y\}}{\max\{x,y\}},~\forall~x,y\in X$ and $t>0,$ is a non-compact, weak $G$-complete fuzzy metric space. 

We now show that $(X,M,\cdot)$ is, \textit{in fact}, Lebesgue. 

Choose a fuzzy pseudo-Cauchy sequence $(x_n)$ of distinct terms in $(X,M,\cdot).$ 

Clearly, $(x_n)$ cannot be eventually in $\left\{\frac{1}{2^n}:n\ge2\right\}.$ For otherwise, there exists $k\in\mathbb N$ such that $M(x_m,x_n,t)\le\frac{1}{2},~\forall~n\ge k,t>0,$ which is a contradiction since $(x_n)$ is pseudo-Cauchy.

So there exists a subsequence $(x_{r_n})$ of $(x_n)$ such that $x_{r_n}\in\left[\frac{1}{2},1\right],~\forall~n\in\mathbb N.$ 

Since $\tau_M$ defines the usual topology of $\mathbb R$ restricted to $X$ \cite{uc}, so $[\frac{1}{2},1]$ is a compact subset of $(X,M,\cdot).$ Consequently, in view of Theorem \ref{uni2}, $(x_{r_n})$ (and hence $(x_n)$) has a cluster point in $X.$

Thus $(X,M,\cdot)$ is Lebesgue.
\end{exm}

\begin{obs}
\normalfont $\mathcal{K}\subsetneq\mathcal{L}.$
\end{obs}

The following result is an immediate consequence of Lemma \ref{uni1}:

\begin{theo}
\normalfont A precompact, weak $G$-complete fuzzy metric space is Lebesgue.
\end{theo}

In view of the last theorem, we have the next corollary from Proposition \ref{prec}, Proposition \ref{pro2} and Proposition \ref{pro1}:

\begin{cor}
\normalfont A precompact, weak $G$-complete metric space is Lebesgue.
\end{cor}


\begin{thebibliography}{99}

\bibitem{2} S. Adhya, A. Deb Ray, On Lebesgue Property for Fuzzy Metric Spaces, TWMS J. of Apl. \& Eng. Math., \textit{to appear}.

\bibitem{gb} G. Beer, Metric Spaces on which Continuous Functions are Uniformly Continuous and Hausdorff Distance, Proc. Amer. Math. Soc. 95.4 (1985) 653-658.

\bibitem{ver1} A. George, P. V. Veeramani, On Some Results in Fuzzy Metric Spaces, Fuzzy Sets and Systems 64 (1994) 395-399.

\bibitem{ver2} A. George, P. V. Veeramani, On Some Results of Analysis for Fuzzy Metric Spaces, Fuzzy Sets and Systems 90 (1997) 365-368.

\bibitem{ogc} M. Grabiec, Fixed Points in Fuzzy Metric Spaces, Fuzzy Sets and Systems 27 (1988) 385-389.

\bibitem{com} V. Gregori, A. López-Crevillén, S. Morillas, A. Sapena, On Convergence in Fuzzy Metric Spaces, Top. Appl. 156 (2009) 3002-3006.

\bibitem{com2} V. Gregori, J. J. Mi$\widetilde{\text{n}}$ana, S. Morillas, Uniform Continuity in Fuzzy Metric Spaces, Rend. Ist. Mat. Univ. Trieste 32 Suppl. 2 (2001) 81-88.

\bibitem{star} V. Gregori, J. J. Mi$\widetilde{\text{n}}$ana, A. Sapena, Completable Fuzzy Metric Spaces, Top. Appl. 225 (2017) 103-111.

\bibitem{gc} V. Gregori, J. J. Mi$\widetilde{\text{n}}$ana, A. Sapena, On Banach Contraction Principles in Fuzzy Metric Spaces, Fixed Point Theory 19 (2018) 235-248.

\bibitem{stat} V. Gregori, S. Romaguera, Characterizing Completable Fuzzy Metric Spaces, Fuzzy Sets and Systems 144 (2004) 411-420.

\bibitem{uni} V. Gregori, S. Romaguera, Some Properties of Fuzzy Metric Spaces, Fuzzy Sets and Systems 115 (2000) 485-489.

\bibitem{uc} V. Gregori, S. Romaguera, A. Sapena, Some Questions in Fuzzy Metric Spaces, Fuzzy Sets and Systems 204 (2012) 71-85.

\bibitem{kmf}  O. Kramosil, J. Michalek, Fuzzy Metric and Statistical Metric Spaces, Kybernetica 11 (1975) 326-
334.

\bibitem{pc} A. Sapena, A Contribution to the Study of Fuzzy Metric Spaces, Applied General Topology 2.1 (2001) 63-75.

\bibitem{tn} B. Schweizer, A. Sklar, Statistical Metric Spaces, Pacific J. Math. 10 (1960) 314-334.

\bibitem{tir} P. Tirado, On Compactness and $G$-Completeness in Fuzzy Metric Spaces, Iranian Journal of Fuzzy Systems 9.4 (2012) 151–158.

\bibitem{will} S. Willard, General Topology. Reading, MA: Addison-Wesley Publishing (1970).
\end{thebibliography}
\end{document}